\documentclass{amsproc}

\newtheorem{theorem}{Theorem}[section]
\newtheorem{lemma}[theorem]{Lemma}
\newtheorem{corollary}[theorem]{Corollary}
\newtheorem{proposition}[theorem]{Proposition}

\theoremstyle{definition}
\newtheorem{definition}[theorem]{Definition}

\theoremstyle{remark}
\newtheorem{remark}[theorem]{Remark}

\numberwithin{equation}{section}

\usepackage{epsfig}

\theoremstyle{plain}
\newtheorem{defiprop}[theorem]{Definition and Proposition}
\newcommand{\C}{{\mathbb{C}}}
\newcommand{\N}{{\mathbb{N}}}
\newcommand{\R}{{\mathbb{R}}}
\newcommand{\Z}{{\mathbb{Z}}}
\newcommand{\rb}{{$b$}} 
\newcommand{\interior}[1]{\overset{\circ}{#1}}

\DeclareMathOperator{\nullfaces}{null}
\DeclareMathOperator{\spec}{spec}

\DeclareMathOperator*{\asym}{\sim} 
\DeclareMathOperator*{\ecup}{\overline{\cup}} 
\DeclareMathOperator*{\bigecup}{\overline{\bigcup}} 
\DeclareMathOperator*{\LIM}{LIM}
\newcommand{\dregint}{-\kern-1.04em\int}
\newcommand{\tregint}{-\kern-0.9em\int}

\begin{document}

\title{Singular Asymptotics Lemma and Push--Forward Theorem}

\author{Daniel Grieser}
\address{Institut f\"ur Mathematik, Humboldt-Universit\"at zu Berlin, Unter den Linden 6, 10099 Berlin, Germany}
\email{grieser@mathematik.hu-berlin.de}
\thanks{DG was supported by the 
Deutsche Forschungsgemeinschaft.}

\author{Michael J.\ Gruber}
\address{Department of Mathematics, M.I.T. 2--167, 77 Massachusetts Avenue, Cambridge, MA 02139--4307, USA}
\email{mjg@math.mit.edu}
\thanks{MG was supported in part by Deutsche Forschungsgemeinschaft, SFB 288, project D6.}

\subjclass{Primary 41A60,58G18; Secondary 30E15, 33A40, 35A27,  35S05, 35S35,  47F05, 47G30, 58G10, 58G11, 58G15}
\date{December 15, 1999 and, in revised form, \today.}

\begin{abstract}
The Singular Asymptotics Lemma by Br\"uning and Seeley and the Push-Forward
 Theorem by Melrose lie at the very heart of their respective approaches to
 singular analysis. We review both and show that they deal with the same
basic problem, giving solutions that emphasize different aspects of it.
This also points to a possible common extension.
\end{abstract}

\maketitle


\section{An Example} \label{MG-sec:Intro}
\begin{figure}[htbp]
\[\epsfig{file=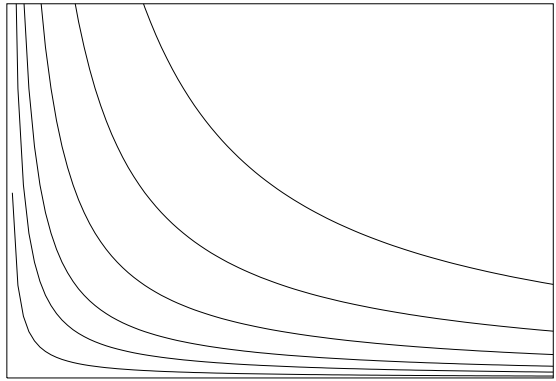,width=3in}\]
\caption{Level lines of $f:(x,y)\mapsto xy$}\label{MG-fig:xy}
\end{figure}

As a simple example, consider 
\begin{equation} f:\R_+^2\rightarrow\R_+, (x,y)\mapsto xy. \label{MG-eq:xy} \end{equation}
The level lines $\{f=c\}$ of $f$ are shown in Figure~\ref{MG-fig:xy}. As $c$ tends to $0$, the corresponding level line approaches the boundary $\{0\}\times\R_+\cup\R_+\times\{0\}$ (where we write $\R_+=[0,\infty)$).

We wish to integrate a compactly supported smooth density
 over the level lines of $f$, i.e.\ to compute the push-forward of $u\,dxdy$
 under $f$, for some $u\in C^\infty_c(\R_+^2)$ 
(that is, $u$ is the restriction to $\R_+^2$ of a compactly supported function
 on $\R^2$): we pair  with $v\in\C^\infty_c(\R_+)$ to get
\begin{align}
 \left\langle f_*(u\,dxdy), v\right\rangle &= \left\langle u\,dxdy,f^*v\right\rangle 
 = \iint_{\R_+^2} u(x,y)v(xy)\,dxdy \notag \\
 &= \int_{\R_+} \left( \int_{\R_+} u\left(x,\frac tx\right)\frac1x\,dx\right) v(t)\,dt, \\
\intertext{so that}
f_*u(t) &= \int_{\R_+} u\left(x,\frac tx\right)\frac1x\,dx \label{MG-eq:f*u}
\end{align}
in distributional sense\footnote{%
Recall that the push-forward of a distribution (under a smooth proper map) always exists as a distribution, and that the push-forward of a smooth function may be non-smooth. On the other hand, the pull-back of a smooth function is always smooth, but not every distribution can be pulled back under smooth proper maps.

In the following we identify densities with functions by using the standard (Lebesgue) measure on $\R^2$ and $\R$, as we did already in \eqref{MG-eq:f*u}.
}. 
$u$ is smooth by assumption and has therefore an asymptotic expansion at the boundary of $\R_+^2$ given by the Taylor series, but $f_*u$ may be non-smooth: Choose a smooth cut-off function $\rho_\varepsilon:\R_+\rightarrow[0,1]$ which is $1$ on $[0,1]$ and $0$ on $[1+\varepsilon,\infty)$.
If we set $u_\varepsilon(x,y)=\rho_\varepsilon(x)\rho_\varepsilon(y)$, then the push-forward satisfies:
\begin{align}
f_*u_\varepsilon(t) &=  \int_{t/(1+\varepsilon)}^{(1+\varepsilon)} \rho_\varepsilon(x)\rho_\varepsilon\left(\frac tx\right)\frac1x\,dx 
  \begin{cases} \leq (1+\varepsilon)\ln(1+\varepsilon) - \ln t \\
   \geq - \ln t \end{cases}
\end{align}
for $t>0$. In other words: Already the innocent constant $1$ (cut off by $\rho_\varepsilon(x)\rho_\varepsilon(y)$) produces a logarithmic singularity, when pushed forward under the seemingly innocent $f$. 
Note that $f_*u_\varepsilon(t)$ tends to $-\ln t$ pointwise as $\varepsilon$ tends to $0$, and also in the distributional sense. Therefore, $f_*u_0(t)=-\ln t$, where $u_0$ is the characteristic function of $[0,1]^2$.

This raises  the following questions:
\begin{itemize}
\item What happens for more general $u$?
\item What happens for more general $f$?
\end{itemize}
The Push-Forward Theorem (PFT) answers these questions by describing the type of
 asymptotic behavior of the push-forward which can occur when $f$ is a
 \rb-fibration and $u$ is a polyhomogeneous conormal distribution on a compact 
manifold with corners.

\bigskip
Another way to look at \eqref{MG-eq:f*u} is to interpret it as integration of asymptotics: Set 
\begin{align} \sigma(x,\zeta)&=u\left(x,\frac1\zeta\right)\frac1\zeta \label{MG-eq:sigma} \\
\intertext{and $z=1/t$, then}
f_*u(t)&=z\int \sigma(x,zx)\,dx. \label{MG-eq:sigmaint}
\end{align}
 Since $u\in C^\infty_c(\R_+^2)$, the asymptotic behavior of $\sigma(x,\zeta)$ for small $x$ and for large $\zeta$ is given by the Taylor expansion of $u(x,y)$ at $x=0$ and at $y=0$. 
One might ask whether this asymptotic expansion can be integrated to get the asymptotic behavior of $f_*u(t)$ for small $t$, i.e.\ to get an expansion of \eqref{MG-eq:sigmaint} for large $z$. 
For $u\equiv u_0$ as above, i.e.\ $\sigma(x,y)=1/y$ on $[0,1]\times[1,\infty)$ and $\sigma(x,y)=0$ otherwise, one recovers the logarithmic singularity, of course.

This raises the following questions:
\begin{itemize}
\item What happens for more general $u$ (or $\sigma$)?
\item What is the structure of the coefficients in the expansion?
\end{itemize}
The Singular Asymptotics Lemma (SAL) answers these questions by computing the coefficients in terms of derivatives and generalized integrals of $\sigma$, for $\sigma$ having quite general asymptotics.

In Section 2 we review the SAL and in Section 3 the PFT. 
In Section 4 we analyze the example above from both the SAL and PFT
point of view. In Section 5 we drop the assumption of compact $y$-support
for $u$ and allow more general asymptotics as $y\to 0$. This is essentially
the generality treated by SAL. We construct the corresponding PFT setting
and show that the 'integrability conditions' in SAL and PFT correspond.
A short summary and some remarks in Section 6 on the different perspectives of SAL and PFT conclude
the paper.


\section{The Singular Asymptotics Lemma}
This lemma was set up in \cite{MG-BruSee:RSA} and used to study singular differential operators defined by the Friedrichs extension of $-\partial_x^2+x^{-2}a$ on $L^2(0,\infty)$, where $a$ is smooth and $a\geq-\tfrac14$, $|a(x)|\leq C(1+x)$. In the sequel this has been used to determine heat trace asymptotics for conical singularities (\cite{MG-BruSee:RESORSO}), which corresponds to allowing $a$ to be a family of unbounded operators. For this, and further references, as well as generalizations to operators of higher orders, we refer to \cite{MG-Les:OFTCSAM}, which we follow closely in this section. We outline the main ideas, but do not reproduce the proofs.

To state SAL's assertion about the coefficients appearing in the expansion we need to define a regularized integral.
We do this for functions having finite asymptotic expansions in the $L^1$-sense:

\begin{definition} \label{MG-defi:Lpq}
Let $p,q\in\R$, $p>0,q>0$. We define ${L}_{p,q}(\R_+^\times)$, where  $\R_+^\times=(0,\infty)$, to be the set of functions $f$ for which there are complex numbers $\alpha_j,\beta_j\in\C$ and (nontrivial) polynomials $p_{\alpha_j},q_{\beta_j}\in\C[\ln x]$   such that
\begin{subequations}\begin{align}
f(x) &= \sum_{j=1}^N x^{\alpha_j}p_{\alpha_j}(\ln x) +x^pf_1(x) \label{MG-eq:pj}\\
 &=\sum_{j=1}^Mx^{\beta_j}q_{\beta_j}(\ln x) +x^{-q}f_2(x) \label{MG-eq:qj}
\end{align}\end{subequations}
for some $f_1\in L^1_{loc}([0,\infty)),f_2\in L^1([1,\infty))$. We assume that the real parts $\Re \alpha_j$ are increasing and $\leq p-1$, and that $\Re \beta_j$ is decreasing and $\geq -q-1$.
Furthermore, we define the set of exponents as
\begin{subequations}\begin{align}
\spec_0(f) &= \{ \alpha_j\mid j=1,\ldots,N \}, \\
\spec_\infty(f) &= \{ \beta_j\mid j=1,\ldots,M \}.
\end{align}\end{subequations}
Finally we assume that all $\alpha_j$ are distinct and set, for arbitrary $\alpha\in\C$, $p_\alpha\equiv p_{\alpha_j}$ if $\alpha=\alpha_j$ and $p_\alpha\equiv0$ if there is no such $\alpha_j$ (analogous for $q_\beta$).
\end{definition}

Now, the shortest way to define the regularized integral of functions in ${L}_{p,q}(\R_+^\times)$ is by using the following limit in the mean:
\begin{definition}
For $f\in{L}_{p,q}(\R_+^\times)$, $p>0,q>0$,  with asymptotics \eqref{MG-eq:pj} and \eqref{MG-eq:qj} we define
\begin{subequations}\begin{align}
\LIM_{x\to0} f(x) &:= p_0(0), \\
\LIM_{x\to\infty} f(x) &:= q_0(0).
\end{align}\end{subequations}
\end{definition}
In other words: the limit in the mean is defined as the constant term in the asymptotic expansion.

\begin{defiprop}
Let $p>0,q>0$, $f\in{L}_{p,q}(\R_+^\times)$. If we define 
\[ F(x):=\int_1^xf(t)\,dt \]
then $F\in{L}_{p+1-\varepsilon,q+1-\varepsilon}(\R_+^\times)$ for every $\varepsilon>0$.

We define the regularized integral $\tregint f$ of $f$ by
\begin{equation}
\dregint_0^\infty f(x)\,dx := \LIM_{x\to\infty}F(x)- \LIM_{x\to0}F(x).
\end{equation}
\end{defiprop}

\begin{remark} \hglue0pt

\noindent\begin{itemize}
\item A more conceptual approach is to use the Mellin transform ${\mathcal M}f(z)=\int_0^\infty x^{z-1}f(x)\,dx$. A priori it is defined for arbitrary $z\in\C$ in case of $f\in C_c^\infty(\R_+^\times)$, and ${\mathcal M}f(1)=\int_0^\infty f$. 

As it turns out, functions $f\in{L}_{p,q}(\R_+^\times)$ still have a well-defined Mellin transform which is meromorphic for $1-p< \Re z<1+q$, and has at most poles of order $\deg p_\alpha +1$ in $-\alpha$ for $\alpha\in\spec_0 (f)$, and poles of order $\deg q_\beta+1$ in $-\beta$ for $\beta\in\spec_\infty (f)$. 
Then it is natural to regularize $\int$ by defining $\tregint f$ to be the zeroth order residuum of ${\mathcal M}f(z)$ at $z=1$. This leads to the same definition as above.
\item An easy calculation shows $\tregint_0^\infty x^\alpha\ln^kx\,dx=0$, since the contributions at $0$ and at $\infty$ just cancel.
\end{itemize}
\end{remark}

{}From the definition, $\tregint$ shares the linearity properties with $\int$. But a fundamental difference appears in the substitution rule:

\begin{lemma}
Let $p>0,q>0$, $f\in{L}_{p,q}(\R_+^\times)$, $t >0$. Set $P_{-1}(\xi):=\int_0^\xi p_{-1}(s)\,ds$ and $Q_{-1}(\xi):=\int_0^\xi q_{-1}(s)\,ds$ for the coefficient polynomials $p_{-1},q_{-1}$ from the asymptotic expansions \eqref{MG-eq:pj} and \eqref{MG-eq:qj}. Then
\begin{equation}
\dregint_0^\infty f(t  x)\,dx=\frac1t \left(\dregint_0^\infty f(x)\,dx + Q_{-1}(\ln t )-P_{-1}(\ln t )  \right).
\end{equation}
\end{lemma}

By applying this lemma to a function $\varphi(tx) f(x)$, where $f$ is as above and $\varphi$ is a Schwartz function (and therefore has asymptotics given by its Taylor series), we arrive at the first SAL--like result:
\begin{proposition}
Let $p>0,q>0$, $f\in{L}_{p,q}(\R_+^\times)$, $\varphi\in{\mathcal S}(\R)$. Then, as $t\to0$,
\begin{equation}\begin{aligned}
\dregint_0^\infty\varphi(tx)f(x)\,dx&= \sum_{j<q} \frac{\varphi^{(j)}(0)}{j!} \dregint_0^\infty x^j f(x)\,dx \cdot t^j \\
&+ \sum_{\substack{\beta\in\spec_\infty (f)\\\Re\beta\geq-q-1}} \dregint_0^\infty\varphi(x)x^\beta q_\beta\left(\ln\frac xt\right)\,dx\cdot t^{-\beta-1} \\
&+ \sum_{\substack{\beta\in\spec_\infty(f)\cap\Z\\-1\geq\beta\geq-q-1}} \frac{\varphi^{(-\beta-1)}(0)}{(-\beta-1)!} Q_\beta(-\ln t)\cdot t^{-\beta-1}+{\mathcal O}(t^q).
\end{aligned}\end{equation}
\end{proposition} 
The poles from $\spec_0(f)$ produce additional terms if we study $\varphi(x)f(\tfrac xt)$ by applying the substitution lemma to the previous proposition:
\begin{corollary}
Let $\varphi,f$ be as above. Then, as $t\to0$,
\begin{equation}\begin{aligned}
\dregint_0^\infty\varphi(x)f\left(\frac xt\right)\,dx&= \sum_{j<q} \frac{\varphi^{(j)}(0)}{j!} \dregint_0^\infty x^j f(x)\,dx \cdot t^{j+1} \\
&+ \sum_{\substack{\beta\in\spec_\infty (f)\\\Re\beta\geq-q-1}} \dregint_0^\infty\varphi(x)x^\beta q_\beta\left(\ln\frac xt\right)\,dx\cdot t^{-\beta} \\
&+ \sum_{\substack{\beta\in\spec_\infty(f)\cap\Z\\-1\geq\beta\geq-q-1}} \frac{\varphi^{(-\beta-1)}(0)}{(-\beta-1)!} Q_\beta(-\ln t)\cdot t^{-\beta} \\
&- \sum_{\substack{\alpha\in\spec_0(f)\cap\Z\\-1\geq\alpha\geq-q-1}} \frac{\varphi^{(-\alpha-1)}(0)}{(-\alpha-1)!} P_\alpha(-\ln t)\cdot t^{-\alpha}+{\mathcal O}(t^{q+1}).
\end{aligned}\end{equation}
\end{corollary}
If we pass now from $\varphi\otimes f$ to more general functions
 $\sigma$ of two variables, then it is easy to derive
\begin{theorem}[Singular Asymptotics Lemma]
Let $C$ be the sector $\{\xi\in\C^\times\mid |\arg\xi|<\pi-\varepsilon\}$ for some $\varepsilon>0$. Assume that $\sigma:\R\times C\to\C$ fulfills the following properties for some $p$:
\begin{enumerate}
\item $\sigma$ is $(p-1)$-times continuously differentiable in the first variable $x$, and $\partial_x^{(p-1)}\sigma(\cdot,\zeta)$ is absolutely continuous on $[0,\infty)$.
\item All derivatives up to order $p$ are analytic\footnote{Even in the analytic case assumption \eqref{MG-eq:fp-1} below is a non-empty condition, since the sector $C$ does not contain $0$.} (resp.\ measurable\footnote{In the measurable case the asymptotics will only hold along the real line $z\to\infty$, not inside the whole sector $z\in C$.}) in the second variable $\zeta$.
\item There are polynomials $p_\alpha\in {\mathcal S}(\R)[\ln\zeta]$ with Schwartz functions as coefficients, and constants $C_{JK}$ and $r$ such that
\begin{align} \left\vert x^J\partial_1^K\left(\sigma(x,\zeta)-\sum_{\Re\alpha>-p-1} \zeta^\alpha p_\alpha(x,\ln\zeta)\right)\right\vert&\leq C_{JK}|\zeta|^{-p-1}|\ln^r\zeta| \label{MG-eq:SALass}\end{align}
for $|\zeta|\geq1$, $0<x\leq|\zeta|$, $K\leq p$, $J\in\Z_+$
\item There is $f_p:(0,1]\to\R$ such that 
one has:
\begin{subequations}\begin{align}
\int_0^1\left\vert\zeta^j\partial_1^j\sigma(0,\zeta)\right\vert\,d\zeta &<\infty ,\qquad j=0,\ldots,p-1, \label{MG-eq:fp-1}\\
\int_0^1\int_0^1s^p\left\vert\partial_1^p\sigma(\theta st,s\xi)\right\vert\,ds\,dt &\leq f_p(\theta),\qquad\text{if } p>0, \label{MG-eq:fp>0}\\
\int_0^1\left\vert\sigma(\theta s,s\xi)\right\vert\,ds&\leq f_0(\theta),\qquad\text{if } p=0, \tag{\ref{MG-eq:fp>0}'} \label{MG-eq:fp=0}
\end{align} \label{MG-eq:SALass2}\end{subequations}
for $0<\theta\leq1$, uniformly for $|\xi|=1$.
\end{enumerate}
Then, as $z\to\infty$ in $C$ (resp.\ for $z\geq0$), one has
\begin{equation}\begin{aligned}
\int_0^\infty\sigma(x,xz)\,dx&= \sum_{j=0}^{p-1}\dregint_0^\infty \frac{\zeta^j}{j!}\partial_1^j\sigma(0,\zeta)\,d\zeta\cdot z^{-j-1} \\
&+\sum_{\Re\alpha>-p-1} \dregint_0^\infty x^\alpha p_\alpha(x,\ln xz)\,dx\cdot z^\alpha \\
&+\sum_{\alpha=-1}^{-p}\frac{z^\alpha}{(-\alpha-1)!}\partial_1^{(-\alpha-1)}P_\alpha(0,\ln z) \\
&+{\mathcal O}\left(z^{-p-1}\ln^{r+1}z\right)+{\mathcal O}\left(z^{-p-1}f_p(|z|^{-1})\right)
\end{aligned}\end{equation}
Here, we defined $P_\alpha(x,\ln\zeta)=\int_0^{\ln\zeta} p_\alpha(x,s)\,ds$; and, for $p\in {\mathcal S}(\R)[\ln\zeta]$, we write $p(x,\ln\zeta)$ where $\ln\zeta$ denotes the polynomial variable.
\end{theorem}

The value of SAL's statement clearly depends on the function $f_p$ appearing in the estimates \eqref{MG-eq:fp-1}--\eqref{MG-eq:fp=0}. In the original work \cite{MG-BruSee:RSA}  $f_p$ was constant. In some situations, even on cones (see \cite{MG-Les:OFTCSAM}) one needs to relax this, by using the following criterion:

\begin{proposition}[supplement to SAL]
Assume that we have for some $\varphi\in C([0,1])$, $\varepsilon>0$, $T>0$, the estimates
\begin{subequations}
\begin{align}
|\sigma(x,\zeta)|&\leq \varphi(x)|\zeta|^{-1}, & 0<|\zeta|\leq1, &0<x\leq1,\\
|\sigma(\theta s,\xi s)| &\leq c\theta^{-T}s^{\varepsilon-1},  &0<\theta\leq 1 ,  &0<s\leq 1
\end{align}
\end{subequations}
uniformly in $|\xi|=1$.
Then
\begin{align}
\int_0^1|\sigma(\theta s,\xi s)|\,ds &= {\mathcal O}(\ln\theta)\text{ as }\theta\to0.
\end{align}
\end{proposition}

This represents the case ``$p=0$'' in SAL and allows to identify terms in the asymptotics up to order $\ln t$.


\section{The Push-Forward Theorem}
Since the \rb-calculus is covered already in another contribution in this
volume \cite{MG-Gri:BBC}
 we will keep this section short compared to the previous one. The main sources for this section are \cite{MG-Mel:POCSL,MG-Mel:CCDMC}.

Let $X$ be a manifold with corners. 
An index set is a non-empty discrete set $K\subset\C\times\N_0$ with
\begin{align}
\begin{gathered}
(\alpha_j,k_j)\in K,\,|(\alpha_j,k_j)|\to\infty\Rightarrow\Re \alpha_j\to\infty,\\
(\alpha,k)\in K,\,p\in\N\Rightarrow (\alpha+p,k)\in K, \\
(\alpha,k)\in K,\,p\in\N_0,\,0\leq p<k\Rightarrow (\alpha,p)\in K.
\end{gathered} \label{MG-eq:indexset}
\end{align}
An index family $\mathcal K$ is an assignment of an index set ${\mathcal K}(H)$ to every boundary hypersurface $H$. 
For a given index family one defines the space of polyhomogeneous conormal functions ${\mathcal A}^{{\mathcal K}}_{\text{phg}}(X)$, and these are the functions $u$ which have an asymptotic expansion
\begin{align}
u\asym \sum_{(\alpha,p)\in{\mathcal K}(H)} \rho^\alpha \ln^p\rho\, v_{\alpha,p}
\end{align}
near every boundary hypersurface $H$, with $\rho$ a defining function of
$H$; the coefficients $v_{\alpha,p}$ are polyhomogeneous conormal on $H$, with
index family  $\mathcal K$ restricted to $H$.
Some care needs to be taken in the definition of $\asym$, see
\eqref{MG-eq:uasymp} for the precise definition in the case needed here.
If $X=\R_+$ then
these $u$ clearly correspond to the functions in ${L}_{\infty,\infty}$ with
complete asymptotics from Definition~\ref{MG-defi:Lpq}, with two boundary
points $0$ and $\infty$. The 
index set at $0$ (resp.\ $\infty$) corresponds to the $\alpha_j$ (resp.\ 
$-\beta_j$) and the exponents of the polynomials $p_{\alpha_j}$ (resp.\
$q_{\beta_j}$), and a set of such exponents generates an index set by 
``completion'' according to~\eqref{MG-eq:indexset}.

Now we define several types of maps between manifolds with corners: 
If $H$ is a boundary face of
codimension $k$ we can choose local coordinates
$x_1,\ldots,x_k,y_1,\ldots,y_{n-k}$ of $X$ such that $x_j\geq0$, and
$H=\{x_j=0\}$ near a point $\bar x$. We do the same `in prime' near $\bar
x'=f(\overline{x})$, if $f$ is a map $X\to Y$ between manifolds with
corners. Then, $f$ is 
called \rb-normal if  these coordinates can be chosen so that
\begin{subequations}
\begin{align}
f^*(x_j') &= \prod_{r\in I_j}x_r^{e(r,j)}\text{ for }j=1,\ldots,k'. \label{MG-eq:b-normal} \\
\intertext{for some disjoint $I_j\subset\{1,\ldots,k\}$. Here, the nonnegative
integers $e(r,j)$ define the exponent matrix of $f$. $f$ is called
 $b$-submersion if the coordinates can be chosen so that}
f^*(y_p') &=y_p\text{ for }p=1,\ldots,n'-k'. \label{MG-eq:b-submersion}
\end{align}
\end{subequations}
$f$ is called  \rb-fibration if 
both \eqref{MG-eq:b-normal} and \eqref{MG-eq:b-submersion} can be achieved simultanuously.
See \cite{MG-Gri:BBC} for a more geometric characterization and motivation
of this definition.

The push-forward of an index family is defined as follows:

\begin{definition}
 If $f:X\to Y$ is a \rb-fibration
and $\mathcal K$ is an index family for $X$ then the index family $f_\#{\mathcal K}$ for $Y$ is defined by
\begin{align}
f_\#{\mathcal K}(H)=\bigecup_{\substack{G\in M_1(X)\\e_f(G,H)\ne0}} \left\{\left(\frac{z}{e_f(G,H)},p\right)\mid (z,p)\in{\mathcal K}(G)\right\}
\end{align} for any $H\in M_1(Y)$.
\end{definition}
The extended union used in the previous definition is the following: Let $K,I$ be index sets, then
\begin{align}
K\ecup I= K\cup I\cup\{(z,p'+p''+1)\mid (z,p')\in K,(z,p'')\in I\}. \label{MG-eq:eunion}
\end{align}

Finally, we define the nullset $\nullfaces(e_f)$ of $f$:
\begin{align}
\nullfaces(e_f)=\{G\in M_1(X)\mid \forall_{H\in M_1(Y)}:e_f(G,H)=0\}
\end{align}

\begin{theorem}[Push-Forward Theorem]
If $f:X\to Y$ is a $b$-fibration between compact manifolds with corners, and if $\mathcal K$ is any index family for $X$ with
\begin{equation}
G\in\nullfaces(e_f) \Rightarrow \Re{\mathcal K}(G)>0 \label{MG-eq:nodelta}
\end{equation}
then push-forward under $f$ gives a mapping
\begin{align}
f_*:{\mathcal A}^{{\mathcal K}}_{\text{phg}}(X;{}^b\Omega)&\rightarrow {\mathcal A}^{f_\#{\mathcal K}}_{\text{phg}}(Y;{}^b\Omega) 
\end{align}
\end{theorem}
(See Section \ref{MG-sec:Equiv} for the meaning of ${}^b\Omega$, i.e.\
$b$-densities.)
There is also a version with finite asymptotics, and for partially polyhomogeneous functions.
Clearly, the theorem extends to non-compact $X$ when pushing forward only
compactly supported densities.

The Push-Forward Theorem shows that additional logarithmic terms at a 
hypersurface $H\subset Y$  may be produced in $f_*\omega$ if there are two 
different hypersurfaces $G_1,G_2\subset X$ with $e_f(G_j,H)\ne0,j=1,2,$ 
(i.e.\ both $G_1$ and $G_2$ are mapped into $H$ by $f$), and if the same 
power occurs in the expansions of $\omega$ at $G_1$ and $G_2$ (at least for generic such $\omega$)
\footnote{Actually, a simple partition of unity argument shows that
additional $\log$-terms are created only if $G_1,G_2$ intersect, as the
$x$- and $y$-axis in Figure \ref{MG-fig:xy} do. It is easy to see that
then the level sets of $f$ look roughly like (a higher-dimensional version
of) those in Figure \ref{MG-fig:xy}. Thus: `$\log$'s are produced by level sets
pushing into a corner.'}. 
This explains the occurence of logarithmic terms when pushing forward densities with purely polynomial asymptotics (e.g.\ smooth densities).


\section{The Example Revisited} \label{MG-sec:ExRev}
Let us now apply the theorems from the previous sections to our introductory example~\eqref{MG-eq:xy}, if we push forward smooth densities. Remember $\sigma(x,y)=u(x,\tfrac1y)\tfrac1y$. Therefore, if $u$ is smooth we have the asymptotic expansion
\begin{align}
\sigma(x,\zeta)\asym_{\zeta\to\infty} \sum_{j\geq0} \frac{1}{j!}\partial_2^ju(x,0)\zeta^{-j-1} \label{MG-eq:taylor}
\end{align}
from the Taylor series. Recast in the language of SAL's assumption \eqref{MG-eq:SALass} this means
\begin{align}
p_\alpha(x,\ln\zeta)=\frac1{(-1-\alpha)!}\partial_2^{(-1-\alpha)}u(x,0),\quad\alpha\in-\N.
\end{align}
Furthermore, $\partial_1^j\sigma(x,\zeta)=\partial_1^ju(x,\tfrac1\zeta)\tfrac1\zeta$, and $P_\alpha(x,\ln\zeta)=p_\alpha(x,0)\ln\zeta$ from the definitions.
Thus, the expansion given by SAL is:
\begin{equation} \label{MG-eq:asymp}
\begin{aligned}
f_*u(t)\cdot t&= \int_0^\infty\sigma(x,xz)\,dx \\
&\asym_{t\to0} \sum_{j\geq0}\dregint_0^\infty \frac{\zeta^{j-1}}{j!}\partial_1^ju(0,\tfrac1\zeta)\,d\zeta\cdot t^{j+1} \\
&+\sum_{j\geq0} \dregint_0^\infty \frac{x^{-1-j}}{j!} \partial_2^ju(x,0)\,dx\cdot t^{j+1} \\
&-\sum_{j\geq0}\frac{t^{j+1}}{(j!)^2}\partial_1^{j}\partial_2^ju(0,0)\ln t.
\end{aligned}\end{equation}
In particular, we reproduce the result for $u\equiv u_0$, the characteristic function of $[0,1]^2$, from the introductory section.
\footnote{
Notice the symmetry of \eqref{MG-eq:asymp} in $x$ and $y$ after
the substitution $y=1/\zeta$ in the first integral.
}

In the language of the \rb-calculus, $X=\R_+^2$ has two boundary hypersurfaces ($x$- and $y$-halfaxis), $Y=\R_+$ has one (the point $0$). 
The map $f$ is a \rb-fibration: Near the corner $(0,0)$ we can use cartesian coordinates, near the hypersurface $\{x=0\}$ (resp.\ $\{y=0\})$ we can use coordinates $xy$ and $y$ (resp.\ $xy$ and $x$) to see that it is a \rb-fibration, and that the exponent map is $1$. 
If $u$ is smooth the corresponding index set is $K=\{(n,0)\mid n\in\N_0\}$ for both hypersurfaces. The union is $K\cup K=K$ of course, but the extended union is
\begin{align}
K\ecup K=\{(n,0)\mid n\in\N_0\}\cup\{(n,1)\mid n\in\N_0\}
\end{align} 
from~\eqref{MG-eq:eunion}. This is the \rb-explanation for the logarithmic terms.


\section{Equivalences} \label{MG-sec:Equiv}

In the previous section we compared SAL and PFT in the following special situation:
The function $u$ on $\R_+^2$ to which PFT is applied is smooth and
has compact support. Correspondingly, the function $\sigma(x,\zeta)$ to
which SAL is applied has the special asymptotics \eqref{MG-eq:taylor} as $\zeta\to\infty$
 and vanishes
whenever $x>C$ or $\zeta<c$, for some positive numbers $c,C$.
The assumptions of SAL allow more general $\sigma$: There is no
restriction on the support, and the asymptotics may be more general.

In this section we construct a setting for PFT which corresponds to this.
In particular, we will see that the integrability conditions 
\eqref{MG-eq:SALass2} of
SAL correspond to condition \eqref{MG-eq:nodelta} of PFT.
For simplicity, we only consider the case of complete asymptotic expansions,
i.e.\ $p=\infty$ in SAL. Also, we assume that $\sigma$ has compact
support in $x$, the 
extension to non-compact support in $x$ (but
Schwartz function behavior as expressed in \eqref{MG-eq:SALass})
being  straightforward.

We need to find a manifold with corners $X$ and a \rb-fibration
$F: X\to \R_+$ such that integration $\int \sigma(x,xz) \,dx$ corresponds
to a pushforward by $F$ and the asymptotics $z\to\infty$ corresponds
to the conormal expansion at $0\in\R_+$. The latter suggests introducing
$t=1/z$ as before. Integration $\int h(x,t)\,dx$ of a function $h$
on $\R^2_+$ is just the push-forward under the map $\pi_2:\R^2_+\to\R_+,
(x,t)\mapsto t$. In our case, $h(x,t) = \sigma(x,x/t)$ behaves badly
at $(0,0)$, and the occurence of $x/t$ here suggests blowing up this
point. Thus, set
\begin{equation}\label{MG-eq:Fdef}
X = [\R^2_+,(0,0)], \quad F = \pi_2 \circ \beta : X\to\R_+
\end{equation}
where $\beta:X\to \R^2_+$ is the blow-down map\footnote{
 $[\R^2_+,(0,0)]$ may be defined simply as $[0,\pi/2]\times \R_+$, with
blow-down map polar coordinates: $\beta(\theta,r)=(r\cos\theta,
r\sin\theta)$. While $r,\theta$ have the advantage of providing global
coordinates on $X$, we prefer to avoid transcendental functions
and rather use local coordinates naturally suggested by the problem,
following R.\ Melrose's tradition.
}.
See Figure~\ref{MG-fig:blowdown} for the level sets of $F$.
\begin{figure}[htbp]
\input{blow.pstex_t}
\caption{The blow up $[\R^2_+,(0,0)]$ of $\R^2_+$}\label{MG-fig:blowdown}
\end{figure}

The coordinate functions $x,t$ can be considered (via pullback by $\beta$)
as functions on $X$, and on the interior $\interior{X}$ we also have the 
functions
\begin{equation}\label{MG-eq:coord}
 y = t/x,\quad \zeta = x/t = 1/y.
\end{equation}
Referring to Figure~\ref{MG-fig:blowdown} for notations
(where $G_1,G_2,G_3$ denote the boundary lines and $A=G_1\cap G_3,
B=G_2\cap G_3$ the corners),
we see that $(x,y)$ provide a local coordinate system for $X$ near $A$
(i.e. that they extend smoothly and with independent differentials
to $X\setminus G_2$ and that
$x,y$ are boundary defining functions for $G_3,G_1$ respectively),
and that $(\zeta,t)$ are coordinates near $B$ (i.e. extend to $X\setminus G_1$
and define $G_2,G_3$ respectively).

We now use \rb-densities on $X$ since there is no `canonical' density
to identify functions with densities and since PFT is formulated in
terms of \rb-densities (rather than densities).
By definition, a \rb-density $u\in {\mathcal F}(X,{}^b\Omega)$
(for some function space ${\mathcal F}$) on $X$ can be written
$$u=u_A(x,y) \frac{dx}x\frac{dy}y$$ near $A$ and
$$u=u_B(\zeta,t)\frac{d\zeta}\zeta \frac{dt}t$$ near $B$ (omitting absolute
value signs for simplicity), with $u_A,u_B\in{\mathcal F}$. Note that
on $\interior{X}$ we have 
$\frac{dx}x\frac{dy}y = \frac{d\zeta}\zeta \frac{dt}t$ from 
\eqref{MG-eq:coord}
so that $u_A$ and $u_B$ simply transform like functions, i.e.
\begin{equation} \label{MG-eq:uAuB}
 u_B(\zeta,t) = u_A(\zeta t,1/\zeta).
\end{equation}
(If one would use just densities, there would be additional factors
occuring here.)

To finish the description of the PFT-setting,
we define index sets for the three boundary faces of $X$:
Let $K_1$ and $K_2$ be any index sets and $K_3=\{(n,0): n=1,2,3,\ldots\}$
be the index set of smooth functions vanishing to first order at $G_3$.
Note that for the map $F$ above we have $(e_F(G_i,H))_{i=1,2,3} = (1,0,1)$
where $H=\{0\}$ is the only boundary face of $\R_+$. Therefore, $F$ is 
a $b$-fibration and
$\nullfaces(e_F) = \{G_2\}$, so we have
$$ \text{ \eqref{MG-eq:nodelta} is satisfied } \Leftrightarrow \Re K_2 > 0.$$
The correspondence between SAL and PFT can then be expressed as follows:

\begin{proposition} \label{MG-prop:SALPFT}
Define $X$ and $F$ as in \eqref{MG-eq:Fdef}.
Let $u$ be a compactly supported b-density on $X$, smooth in the interior
$\interior{X}$. 
Define $\sigma:(\R_+^\times)^2\to\R_+$ by
\begin{equation} \label{MG-eq:sigmadef}
 \sigma (x,\zeta) = \frac1x u_A\left(x,\frac1\zeta\right). 
\end{equation}
 
Then 
\begin{equation}\label{MG-eq:formel}
F_*u = \left(\int_0^\infty \sigma(x,\frac{x}t) \,dx\right)\, \frac{dt}t,
\end{equation}
and we have the following equivalences:
\begin{enumerate}
\item[(A)] $\sigma$ is smooth in $\{x\geq0,\zeta>0\}$ $\iff$
$u$ is polyhomogeneous at $\interior{G_3}$, with index set $K_3$
(i.e. $u=xv$ with $v$ smooth up to $\interior{G_3}$), 
\item[(B)] 
$\sigma$ has complete asymptotics, i.e.
\begin{equation} \label{MG-eq:sigmaasymp}
\left|\partial_1^K \left(\sigma(x,\zeta) - \sum_{(\alpha,k)\in K_1,
\Re \alpha < N} \zeta^{-\alpha} \ln^k \zeta\, \sigma_{\alpha k}(x)\right)\right| 
\leq
C_{NK} \zeta^{-N} 
\end{equation}
 for all $N,K$, $x>0, \zeta>1$, with
$\sigma_{\alpha k}\in C_0^\infty(\R_+)$
$$\iff$$
$$ u \text{ is polyhomogeneous, with index sets $K_1,K_3$,
on $X\cap\{y<1\}$ (i.e.\ `near $A$')}, $$
\item[(C)]
If $u\in {\mathcal A}^{(K_1,K_2,K_3)}(X)$ then 
\begin{equation} \label{MG-eq:condition4}
 \int_0^1 \zeta^p |\partial_1^p \sigma(\zeta t,\zeta)|\,d\zeta \leq C_p 
\text{ for all $p$ and $t\in [0,1]$}
\end{equation}
$$\iff$$
condition \eqref{MG-eq:nodelta} of the PFT is satisfied for $u$.
\end{enumerate}
\end{proposition}
\begin{remark}
Roughly speaking, (B) and (C) say that conditions 3. and 4. of SAL correspond
to conditions on the behavior of $u$ near $A$ and $B$. More precisely,
the first halves of (A), (B) and (C) are analogues of 1.+2.,
3. and 4. in SAL for the case $p=\infty$, when $\sigma$ has compact
$x$-support and $f_p=$const, for the case of real $\zeta$ (but SAL 4
is a little weaker than \eqref{MG-eq:condition4} plus polyhomogeneity,
see below).

In (C) we assume that $K_2$ is the `minimal' index set for $u$,
i.e. that the coefficients in the corresponding expansion are non-trivial.
\end{remark}

\begin{proof}[Proof of \ref{MG-prop:SALPFT}]
\eqref{MG-eq:formel} is just \eqref{MG-eq:f*u} in the present notation.

(A) This is immediate from \eqref{MG-eq:sigmadef} since $x$ defines
$G_3$ outside $B$. Note that
the boundary points $A,B$ of $G_3$ correspond to $\zeta=\infty$
and $\zeta=0$, respectively.

(B)  Polyhomogeneity of $u$ near $A$ 
means asymptotic expansions at the interiors
of $G_1, G_3$ (i.e. for $y\to 0$ with $x$ bounded away from 0 and $\infty$,
and vice versa), plus uniformity (`joint asymptotics') as both $x,y$ approach
zero. By definition, this means that, for all $N$,
\begin{equation} \label{MG-eq:uasymp}
u_A(x,y) - \sum_{(\alpha,k)\in K_1, \Re \alpha<N} y^\alpha\ln^k y\, u_{\alpha k}(x)
\in \dot{C}^N_{G_1}(X),
\end{equation}
the space of $N$ times differentiable functions (in $x\geq0,y\geq 0$)
 vanishing to order $N$
at $G_1 = \{y=0\}$, with $u_{\alpha k}\in {\mathcal A}^{K_3}(G_3)$.
Since $y=1/\zeta$, it is easily seen that this is equivalent
to \eqref{MG-eq:sigmaasymp}.

(C) From \eqref{MG-eq:uAuB} and \eqref{MG-eq:sigmadef}
we have 
$$ u_B(\zeta,t) = \zeta t\, \sigma(\zeta t,\zeta).$$
Setting $v_B(\zeta,t) = u_B(\zeta,t)/t$ and taking $t$-derivatives, we get
$$ \partial_2^p v_B(\zeta,t) = \zeta^{p+1} \partial _1^p \sigma(\zeta t,\zeta),$$
so \eqref{MG-eq:condition4} is equivalent to
\begin{equation} \tag{\ref{MG-eq:condition4}'} \label{MG-eq:condition4'}
  \int_0^1  |\partial^p_2 v_B(\zeta,t)|\, \frac{d\zeta}\zeta \leq C_p.
\end{equation}
Now polyhomogeneity of $u$ near $B$ with the index sets $K_2,K_3$ implies
that $\partial_2^p v_B \in {\mathcal A}^{K_2}(G_2)$ for each $p$ and each $t\geq0$.
\eqref{MG-eq:condition4'} holds precisely when the
leading term in the asymptotics $\zeta\to0$ of $\partial_2^p v_B$
has a positive $\zeta$-power, i.e. when $\Re K_2>0$.
\end{proof}

\begin{remark}
 In the proof of (C) very little of the polyhomogeneity at $G_2$
was used. Condition 4 of SAL indicates the precise sense in which this
condition in PFT could be weakened.
\end{remark}
\begin{remark}
 The setup of Sections \ref{MG-sec:Intro} and \ref{MG-sec:ExRev}
 corresponds to the special case
where $u$ is supported in $X\setminus G_2$ and smooth.
\end{remark}


\section{Conclusion}
The basic problem behind both the Singular Asymptotics Lemma and the
Push-Forward Theorem is to determine the asymptotics of integrals like
\eqref{MG-eq:f*u} as $t\to 0$ (or equivalently \eqref{MG-eq:sigmaint}
 as $z\to\infty$), under certain
assumptions on the asymptotics and support of $u$ (resp. $\sigma$).

In the PFT these integrals are interpreted as the push-forward of $u$ under
the map $f$ in \eqref{MG-eq:xy} (resp. $F$ in \eqref{MG-eq:Fdef} when the $y$-support
of $u$ is not compact). The singular nature of this push-forward operation
-- evident from the creation of additional $\log$-terms in the asymptotics
of $f_*u$ -- stems from the fact that $f$ (resp. $F$) is not a fibration
($df$ is not surjective at $(0,0)$, $dF$ at $A$).
The PFT shows that for a certain general class of maps $f$ (\rb-fibrations)
between manifolds with corners, which fail to be fibrations in a controlled
way, the asymptotic type of $f_*u$ can be predicted (or rather
`estimated from above') from that of $u$. The additional $\log$-terms are
predicted from the geometry of the map $f$.

The SAL deals with the integral $\int_0^\infty \sigma(x,zx)\; dx$
directly and therefore gives a more explicit result; 
the singular nature of this integral is reflected in the existence of two
'regimes' ($x>Cz^{-1}$, $x<Cz^{-1}$) in which the integrand behaves in
different ways. 
\footnote{
As pointed out in \cite{MG-Gri:BBC}, `regimes' correspond to minimal
faces of $X$, i.e.\ to the points $A$, $B$. Under the more restrictive
support conditions of Sections \ref{MG-sec:Intro} and \ref{MG-sec:ExRev}
the regime $x<Cz^{-1}$ is empty, so only one corner is needed.
}

We saw in Sections 4 and 5 that the qualitative statements of SAL may be
derived by applying PFT to a special case. 
On the other hand, SAL makes somewhat weaker regularity assumptions and
gives explicit formulas for the coefficients of the resulting expansion;
in particular, this allows to check whether the $\log$-terms (whose
{\em possible} existence is predicted by PFT) actually appear, i.e.
have non-zero coefficients.

\ifx\undefined\allcaps\def\allcaps#1{#1}\fi
\providecommand{\bysame}{\leavevmode\hbox to3em{\hrulefill}\thinspace}


\end{document}